\documentclass[ejs,preprint]{imsart}
\RequirePackage{amsmath,amsfonts,amssymb}
\RequirePackage{graphicx}
\RequirePackage[colorlinks,citecolor=blue,urlcolor=blue]{hyperref}

\startlocaldefs
\def\Cov{\mathrm{Cov}}
\def\Var{\mathrm{Var}}
\newcounter{lem}

\newcounter{cor}
\newcounter{pro}
\newcounter{def}
\newtheorem{theorem}{Theorem}

\newtheorem{corollary}[cor]{Corollary}

\newtheorem{definition}[def]{Definition}

\newtheorem{lemma}[lem]{Lemma}

\newtheorem{proposition}[pro]{Proposition}

\newenvironment{proof}[1][Proof]{\textbf{#1.} }{\ \rule{0.5em}{0.5em}}
\endlocaldefs

\begin{document}

\begin{frontmatter}
\title{A functional Breuer--Major theorem under finite variance, with
application to the $R/S$ test for long memory}
\runtitle{A functional Breuer--Major theorem and the $R/S$ test}
\runauthor{S. Ben Hariz}

\begin{aug}
\author[A]{\fnms{Samir}~\snm{Ben Hariz}\ead[label=e1]{Samir.Ben\_Hariz@univ-lemans.fr}\orcid{0000-0002-1089-0671}}
\address[A]{Laboratoire Manceau de Math\'ematiques, Le Mans Universit\'e,
Avenue Olivier Messiaen, 72085 Le Mans Cedex 09, France\printead[presep={,\ }]{e1}}
\end{aug}

\begin{abstract}
The Breuer--Major theorem (Breuer and Major, 1983) states that, for a
function $G$ of Hermite rank $m$ such that $\mathbb{E}(G^{2}(X))<\infty $,
and a stationary Gaussian sequence whose correlation function $r$ satisfies
$\sum_{i}|r(i)|^{m}<\infty $, the finite dimensional distributions of the
normalized partial sum process of $G(X_{i})$ converge to those of a Brownian
motion. Whether the functional version of this result holds under these sole
conditions has remained open for a long time, the tightness being known only
under additional moment assumptions. In this paper we prove that the
functional central limit theorem holds under finite variance only. More
generally, we give a high order expansion of the partial sum process which
covers both short and long-range dependence: the long-range dependent
components converge to Hermite processes, while the remainder converges to a
Brownian motion.

As an application, we discuss the question of testing long-memory in
observed data using the $R/S$ statistics. We derive the limiting
distributions of the classical and modified $R/S$ statistics, prove the
consistency of the test and propose a data-driven choice of the truncation
parameter. A Monte Carlo study compares the performances of the classical and
modified $R/S$ statistics for both short-memory and long-memory sequences.
\end{abstract}

\begin{keyword}[class=MSC]
\kwd[Primary ]{60F17}
\kwd{60G15}
\kwd[; secondary ]{60G18}
\kwd{62G10}
\kwd{62M10}
\end{keyword}

\begin{keyword}
\kwd{Breuer--Major theorem}
\kwd{functional central limit theorem}
\kwd{Hermite rank}
\kwd{long-range dependence}
\kwd{rescaled range statistic}
\kwd{testing for long memory}
\end{keyword}

\end{frontmatter}

\section{Introduction}

Let $\left( X_{n}\right) _{n\in \mathbb{N}}$ be a stationary sequence
defined on a probability space, and let $G$ be some Borel function. The partial sum
process is defined as:
\begin{equation*}
S_{n}(t,G)=\sum_{i=1}^{[nt]}G(X_{i}),\qquad t\in \lbrack 0,1]
\end{equation*}
and we are interested in the limiting behavior of this process suitably
normalized. When the observed data are i.i.d. with finite second moment, it is
well known that $n^{-1/2}S_{n}(t,G-E(G))$ converges to the Brownian motion,
in the space $\mathcal{D}([0,1])$, endowed with the uniform topology. This
is commonly called Donsker's theorem \cite{don51}. However in many cases, the
assumption of i.i.d. observations may not be valid, especially for applications
in finance, econometrics, hydrology, etc. Several works are concerned with the
generalization of Donsker's theorem to different classes of weakly
dependent processes. Our aim here is to extend these results to long-memory
sequences. Namely, we derive a high order expansion of the Donsker line and
prove that the remainder term converges to a Brownian motion.

Before proceeding further, we summarize some ideas about measuring
dependence in time series analysis.  In general, there are two ways of
measuring dependence :

- With mixing coefficients: this involves the $\sigma $-algebras generated
by the past and the future of the sequence. Let's recall two of the most
commonly used mixing coefficients. Suppose that $(\Omega ,\mathcal{K},%
\mathcal{P})$ is a probability space. For any two $\sigma $-fields $\mathcal{%
A}$ and $\mathcal{B}$ of $\mathcal{K}$, we define the $\alpha $-mixing
coefficient by
\begin{equation*}
\alpha (\mathcal{A},\mathcal{B})=\sup |\mathbb{P}(A\cap B)-\mathbb{P}(A).%
\mathbb{P}(B)|,\quad A\in \mathcal{A},\quad B\in \mathcal{B}
\end{equation*}
and the $\rho $-mixing coefficients by
\begin{equation*}
\rho (\mathcal{A},\mathcal{B})=\sup |\mathrm{corr}(f,g)|,\quad f\in \mathbb{L%
}^{2}(\mathcal{A}),\quad g\in \mathbb{L}^{2}(\mathcal{B}).
\end{equation*}
If $(X_{k})$ is a stationary sequence of random variables, we define
\begin{equation*}
c(n)=\sup_{k}c(\mathcal{F}_{-\infty }^{k},\mathcal{F}_{k+n}^{\infty })
\end{equation*}
where $c=\alpha ,\rho $ and $\mathcal{F}_{n}^{m}$ is the $\sigma $-field
generated by $(X_{i},n\leq i\leq m)$.

- With covariance: this includes essentially associated sequences,
functional of Gaussian and linear processes. For further developments on the subject of dependence and mixing we refer the reader to the survey paper \cite{bra05}.

For the FCLT under strong mixing conditions we refer to the  work by
Merlev\`ede and Peligrad \cite{mer00}. For associated sequences, an optimal
result was obtained by Newman and Wright \cite{new81}. Later in 1993, Birkel
\ \cite{bir93} gave a similar result under a weaker condition on dependence
by introducing the so-called LPQD random variables, however he required a
stronger condition on the moment of $G$. A survey on the FCLT can be found in \cite{Mer06}, see also the monograph \cite{mpu19}. We mention that all these works require at least
the summability of the covariance sequence.

A natural question arises: what will happen if these conditions of weak
dependence are not satisfied. This question is important since many observed
data showed correlations which decay too slow to satisfy the weak \
dependence assumptions. These phenomena are referred to as long-range
dependence or long memory (see \cite{ber94} for more discussions
about long-range dependence, and the monographs \cite{gks12,bfgk13,pt17} for
more recent accounts).

Two families of models are of great importance in the family of long-memory
processes and are widely studied in the past decades.  The first family of
models is moving average or linear processes: They are studied among others
by Giraitis, Surgailis \cite{gir02,gir90}, Robinson \cite{rob03}, Ho and Hsing \cite{hh97}. This family is a natural
generalization of ARMA models which are very popular in time series
analysis.

The second family is functions of Gaussian process which are thoroughly
studied by Arcones \cite{arc94}, Breuer and Major \cite{bre83}, Taqqu \cite{taq79}, Dobrushin and Major \cite{dob79}, and
many others; see also \cite{np12} for a modern approach based on Malliavin
calculus and Stein's method. This class of models has two important features: firstly, it is
a natural generalization of Gaussian processes and i.i.d sequences;
secondly, it can model both dependence in the series and heavy tails in the
marginal distributions. Moreover, application of non-random transformations
on the data may change the behavior of decay in the correlation  function
and this can explain one of the well-known stylized facts of financial data,
namely, a slower decay of the correlations for the absolute value or the
square of the data.

For this second family, the celebrated central limit theorem of Breuer and
Major \cite{bre83} asserts that, if $G$ has Hermite rank $m$,
$\mathbb{E}(G^{2}(X))<\infty $ and $\sum_{i}|r(i)|^{m}<\infty $, then the
finite dimensional distributions of $n^{-1/2}S_{n}(\cdot ,G-\mathbb{E}(G))$
converge to those of a Brownian motion. Surprisingly, the functional version
of this result, that is the tightness of the Donsker line, received much less
attention. Sufficient conditions were given by Chambers and Slud \cite{cs89}
and by Ben Hariz \cite{ben02} for continuous time processes, but these
conditions are not easy to check in practice. Nourdin and Nualart \cite{nn20}
proved the tightness under the additional moment condition
$\mathbb{E}|G(X)|^{p}<\infty $ for some $p>2$, see also Campese, Nourdin and
Nualart \cite{cnn20} for the continuous time case. Very recently, the finite
variance case has been considered by Mansanarez, Poly and Zheng
\cite{mpz26}, under additional assumptions on the covariance function (for
instance $\sum_{i}|r(i)|<\infty $), and by Milesis and Zheng \cite{mz26} by
means of an $\mathbb{L}^{2}$ maximal inequality.

In this paper, we will focus on the question of weak convergence of the
partial sum process for the second family of models. More precisely, let $%
\{X_{n}:n\in \mathbb{N}\}$ be a real valued, Gaussian, stationary sequence
with $\mathbb{E}(X_{0})=0$, $\mathbb{E}(X_{0}^{2})=1$ and a covariance
function $\mathbb{E}(X_{m}X_{n})=r(|m-n|)$. Let $G$ be a real function
satisfying $\mathbb{E}(G^{2}(X))<\infty ,$ where $X$ is a standard normal
variable, we will give an asymptotic expansion for the partial sum
processes,
\begin{equation*}
S_{n}(t,G)=\sum_{i=1}^{[nt]}G(X_{i}),\quad \quad t\in \lbrack 0,1].
\end{equation*}
We show, by a short argument, that the functional central limit theorem holds
under the sole conditions $\mathbb{E}(G^{2}(X))<\infty $ and
$\sum_{i}|r(i)|^{m}<\infty $. Our proof combines a truncation argument with
the maximal inequality established in \cite{bb25}. Moreover, our expansion
handles simultaneously short and long-range dependence: it separates the
long-range dependent components, which converge to Hermite processes, from a
remainder which converges to a Brownian motion. It is the functional
counterpart of the almost sure expansion obtained in \cite{bb25}.

We believe that this development has a potential to give a better insight on
the controversial issue whether financial data has a long memory. This
question was studied by many authors,  see for example  Lo \cite{lo91}, \cite{gir03}, \cite{ban05} and the references therein.
A detailed empirical study is however beyond the scope of this paper. Indeed,
the presence of long memory in financial data, in particular in volatility,
is still debated, and some authors even argue for antipersistence (rough
volatility, see Gatheral et al. \cite{gjr18}), which would require a
two-sided version of the test.

The rest of the paper is structured as follows. In Section 2, we present our
main results. In Section 3, we show how to apply our results to test
long-range dependence in data. In Section 4, we perform Monte Carlo
simulations for both short-range and long-range sequences to study the
performance of the test with the R/S statistic on finite samples. Section 5
is devoted to proofs and technical details. The main tools used in the
proofs are the expansion of functions in Hermite polynomials and moments
inequalities.

\section{Statement of main results}

Before stating our main result we need to introduce some notations. Let $%
\phi $ denotes the normal density and $H_{k}$ the (monic) Hermite polynomial of
degree $k$, defined by $H_{k}\equiv \left( -1\right) ^{k}\phi ^{(k)}/\phi $, so that
$\mathbb{E}\left[ H_{k}(X_{0})H_{l}(X_{i})\right] =k!\,r^{k}(i)\mathbf{1}_{\{k=l\}}$. If $G$ is in $\mathbb{L}^{2}\left( \mathbb{R} ,\phi
\right) ,$ then it can be expanded as:
\begin{equation}
G(x)=\sum_{k=0}^{\infty }{\frac{c_{k}}{{k!}}}H_{k}(x),\qquad c_{k}=\mathbb{E}%
\left[ G\left( X\right) H_{k}(X)\right]  \label{ser}
\end{equation}
where the series (\ref{ser}) is convergent in $\mathbb{L}^{2}$ sense. Here $%
X $ has a density $\phi.$ Let
\begin{equation*}
m\equiv \inf \{k>0:\ \ c_{k}\neq 0\},
\end{equation*}
$m$ is called the \emph{Hermite rank} of $G$. Let also $d_{n,m}^{2}\equiv
\Var\left( S_{n}(1,H_{m})\right).$ When $r(n)\sim n^{-\alpha }L(n),$ where $%
\alpha >0,$ $\ L$ is a slowly varying function and $m\alpha <1,$ then
\begin{equation}
d_{n,m}\equiv \left(\Var\left( S_{n}(1,H_{m})\right) \right) ^{1/2}\sim
\footnote{%
For two sequences we write $u_{n}\sim v_{n}$ if $\lim u_n/v_{n}=1.$}\left(
\frac{2m!}{\left( 1-m\alpha \right) \left( 2-m\alpha \right) }n^{2-m\alpha
}L^{m}(n)\right) ^{1/2}.  \label{dnm}
\end{equation}
Let $\Longrightarrow $ stands for convergence in law in $\mathcal{D}\left(\left[ 0,1\right] \right)$ endowed with the uniform topology. Now we are ready to state our main result.
\begin{theorem}
\label{mth} Let $G$ be a function with Hermite rank $m$ such
that: $\mathbb{E}\left[ G^{2}\left( X\right) \right] <\infty .$ Assume that
there exists an integer $m^{\ast}\geq m$ such that $r^{m^{\ast }}\in \mathbb{L}^{1}$, i.e.
$\sum_{i\in \mathbb{Z}}|r(i)|^{m^{\ast }}<\infty $, and let
\begin{equation}
\sigma ^{2}\equiv \sum_{k=m^{\ast }}^{\infty }{\frac{{c_{k}^{2}}}{{k!}}}%
\sum_{i\in \mathbb{Z}}r^{k}(|i|)\ \geq 0.
\end{equation}
Let $G^{m^{\ast }}$ be defined as
\begin{equation}
G^{m^{\ast }}(x)\equiv \sum_{k=m^{\ast }}^{\infty }{\frac{c_{k}}{{k!}}}%
H_{k}(x).  \label{gms}
\end{equation}
Then
\begin{equation*}
S_{n}(t,G)=\left[ nt\right] \mathbb{E}\left[ G\left( X\right) \right]
+\sum_{k=m}^{m^{\ast }-1}{\frac{{c_{k}}}{{k!}}}S_{n}(t,H_{k})+S_{n}(t,G^{m^{%
\ast }})
\end{equation*}
where
\begin{equation}
n^{-1/2}S_{n}(\cdot ,G^{m^{\ast }})\Longrightarrow \sigma W(\cdot ),
\label{res1}
\end{equation}
in$\,\,\,\mathcal{D}\left( \left[ 0,1\right] \right) ,$ endowed with the
uniform topology and $W(\cdot )$ is the standard Brownian motion (if
$\sigma =0$, the limit is degenerate).

In particular, we have \ the following:

\begin{enumerate}
\item  If $r^{m}\in \mathbb{L}^{1},$ then
\begin{equation}
n^{-1/2}S_{n}(\cdot ,G-\mathbb{E}\left[ G(X)\right] )\Longrightarrow \sigma (G)W(\cdot )  \label{res2}
\end{equation}
where $\sigma ^{2}(G)=\sum_{i\in \mathbb{Z}}\mathrm{Cov}\left(
G(X_{0}),G(X_{|i|})\right) $.

\item  If $r(n)\sim n^{-\alpha }L(n)$ and $m\alpha <1,$ then
\begin{equation}
d_{n,m}^{-1}S_{n}(\cdot ,G-\mathbb{E}\left[ G(X)\right] )\Longrightarrow \frac{{c_{m}}}{m{!}}Y_{m}(\cdot )
\label{res3}
\end{equation}
\ in$\,\,\,\mathcal{D}(\left[ 0,1\right] ),$ endowed with the uniform
topology. Here $Y_{m}(\cdot )$ is the Hermite process of order $m$ (see \cite
{dob79} for a definition).
\end{enumerate}
\end{theorem}

\noindent \textbf{Remark.} The decomposition in Theorem \ref{mth} was used in
\cite[Theorem 3.3]{bb25} to derive almost sure convergence rates
(Marcinkiewicz--Zygmund laws) for each of its components, under possibly
infinite variance. Theorem \ref{mth} is its functional (weak convergence)
counterpart in the finite variance case. The key tool for tightness is the
maximal inequality of \cite[Theorem 2.3]{bb25}, see Lemma \ref{lem1} below.

\noindent \textbf{Comments.} Firstly, we note that our assumptions on the
moment and the correlation function are natural. The finiteness of the second
moment is required for the CLT to hold with the $n^{-1/2}$ normalization even
in the independent setting, thus this condition can not be removed. The
summability condition $\sum_{i}|r(i)|^{m}<\infty $ is the condition of the
Breuer--Major theorem \cite{bre83}, under which the finite dimensional
distributions converge to those of a Brownian motion. Theorem \ref{mth} shows
that, in this discrete time setting, no additional condition is needed for
the tightness. In particular, the extra integrability condition
$\mathbb{E}|G(X)|^{p}<\infty $, $p>2$, used in \cite{nn20} can be dropped.
The summability condition can not be dropped in general: as shown by part 2
of the theorem, when $m\alpha <1$ the normalization and the limiting process
are different. In the boundary case $m\alpha =1$, a Gaussian limit may still
hold, but with a different normalization, and this case is not covered by
Theorem \ref{mth}.

Secondly, we point out that for the finite variance case we showed that the
limiting process depends on the decay of the correlation function. An
interesting question rises: What is the limiting process and the
normalization factor in the infinite variance case, for example, for the
case of symmetric $\alpha $-stable distribution. Some results in this
direction are given in \cite{bb25}, namely almost sure convergence rates
and the convergence to the Hermite process when
$\mathbb{E}|G(X)|^{p}<\infty $ with $(2-m\alpha )\min (2,p)>2$. We believe
that such studies will be helpful for a better understanding of the relation
between heavy tails and long-range dependence, since both factors cause
convergence rate to be slower than the finite variance and\/ or the weak
dependent setting.

Thirdly, we mention that the proof of the functional central limit theorem
(FCLT) relies essentially on a maximal inequality combined with a truncation
argument. We believe that the same scheme can be adapted to other dependence
structures for which suitable maximal inequalities are available, for
instance $\rho $-mixing sequences (see the maximal inequalities of Shao
\cite{sha95} and the CLT of Peligrad \cite{pel87}).

\section{Application to the R/S statistics}

In economics and finance, long-range dependence in time series is an
important and active research topic due to implications of the presence of
long-memory on many paradigms used in \ financial and economics models. For
this reason, the question of how to test long-memory in financial data,
specially in the case of stock returns, has often become a source of major
controversies, see for example \cite{lo91}, \cite{taq99} and \cite{taq99b}.
For an overview on long memory we refer to \cite{ban05} and the references
therein. We also mention the work of  Giraitis et al. \cite{gir03} who considered Rescaled variance and other  tests for long
memory in volatility and levels of returns.

One of the oldest and best-known methods is the so-called $R/S$ statistics.
This method was first introduced by Hurst in 1951 \cite{hur51} while
studying the levels of the Nil river. It was used later by Mandelbrot and
Wallis to analyze long-memory in geophysical data. For a stationary sequence
$\left( Y_{i}\right) _{1\leq i\leq n}$ the rescaled range statistics is
defined as follows:
\begin{equation*}
T\equiv R/S,\quad R\equiv \max_{1\leq k\leq n}\left(
\sum_{i=1}^{k}\left( Y_{i}-\bar{Y}_{n}\right) \right) -\min_{1\leq
k\leq n}\left( \sum_{i=1}^{k}\left( Y_{i}-\bar{Y}_{n}\right) \right)
\end{equation*}
and
\begin{equation*}
S^{2}\equiv S_{n}^{2}=\frac{1}{n}\sum_{i=1}^{n}\left( Y_{i}-\bar{Y}%
_{n}\right) ^{2},
\end{equation*}
where $\bar{Y}_{n}=n^{-1}$ $\sum_{i=1}^{n}Y_{i}$ is the empirical mean and $%
S $ is the standard deviation of the sample. This statistics is used both in
testing and estimating the dependence intensity in a given sequence. In
1991, Lo \cite{lo91} modified the $R/S$ statistics to take into account,
that the limiting variance in the dependent case, involves autocorrelations.
Namely, he replaced the $S^{2}$ by a consistent estimator of $\sigma ^{2}(G),$%
\begin{equation}
S_{n}^{2}(q)\equiv \frac{V_{n}(q)}{q+1}\equiv \frac{1}{n}\sum_{i=1}^{n}%
\left( Y_{i}-\bar{Y}\right) ^{2}+2\sum_{k=1}^{q}\left( 1-\frac{k}{q+1}%
\right) \Cov_{n}(k), \label{snq}
\end{equation}
where $\Cov_{n}(k)\equiv \frac{1}{n}\sum_{i=1}^{n-k}\left( Y_{i}-\bar{Y}%
\right) \left( Y_{i+k}-\bar{Y}\right) .$ Actually $V_{n}(q)$ is the
empirical estimator of \
\begin{equation*}
V(q)\equiv \Var\left( \sum_{i=1}^{q+1}Y_{i}\right) ,
\end{equation*}
the variance of a partial sum of size $q+1,$ i.e the variance of $%
\sum_{i=1}^{q+1}Y_{i}$. Lo derives the limiting behavior of the modified $%
R/S $ for the short-memory case under some weak dependence assumptions
involving mixing coefficients. He concluded that there was no evidence of
long-range dependence in stock prices. The same question was then considered
in several works by Taqqu \textit{et al.} \cite{taq99} and \cite{taq99b}.
These works were essentially based on empirical and simulations studies.
They inferred that the Lo's modified statistics had a strong preference
towards the short-range hypothesis and demonstrated graphically that the
data used in \cite{lo91} had the property of long-memory.

Although the $R/S$ method is widely used and presents some robustness
against non Gaussian marginal distributions with possibly infinite variance,
(see \cite{avr00}), little is known about its asymptotic distribution. It
was reported in \cite{taq99} that ``The less attractive features of the
classical R/S-analysis are its sensitivity to the presence of explicit
short-range dependence structures, its bias, and a lack of a distribution
theory for the underlying statistics''. These are the major drawbacks of
this statistics, and stand as obstacles in the way of using the classical $%
R/S$ as a rigorous statistical inference method. In this section, we will
derive the limiting law of the $T\equiv R/S$ statistics for both
short-memory and long-memory cases. We begin by defining long-range
dependence for our model.

\begin{definition}
Let $\left( Y_{i}\right) =\left( G(X_{i})\right) $ where $\left(
X_{i}\right) $ is a stationary Gaussian sequence and $m$ is the Hermite rank
of $G$.

\begin{itemize}
\item[\textrm{{H0:}}]  The sequence $\left( Y_{i}\right) $ is said to be
short-range dependent (SRD) if $r^{m}\in \mathbb{L}^{1}.$

\item[\textrm{{H1:}}]  The sequence \ $\left( Y_{i}\right) $ is said to be
long-range dependent (LRD) if $r(n)\sim n^{-\alpha }L(n)$ where $\ L$ is a
slowly varying function and $m\alpha <1.$
\end{itemize}
\end{definition}

Roughly speaking, \ the SRD definition is equivalent to $\sum_{i\in \mathbb{Z%
}}\left| \Cov(Y_{0},Y_{i})\right| <\infty ,$ while LRD \ means that $%
\sum_{i\in \mathbb{Z}}\left| \Cov(Y_{0},Y_{i})\right| =\infty .$ Actually,
the SRD case may include the case when $\sum_{i=1}^{n}r^{m}\left( i\right) $
is slowly varying function of $n$. The LRD case may include more general
correlation functions, however the limiting distributions depend on the
decay rate of the covariance function. Note that the boundary case
$r(n)\sim n^{-1/m}L(n)$ (i.e. $m\alpha =1$) is covered neither by H0 nor by
H1 when $\sum_{i}|r(i)|^{m}=\infty $: the partial sums are then of order
$\left( n\sum_{|i|\leq n}|r(i)|^{m}\right) ^{1/2}$, e.g. $\sqrt{n\log n}$
when $L\equiv 1$.

\begin{corollary}
\label{rss} Let $Y=G(X)$ where $X$ is a stationary Gaussian sequence. For
any process $Z$, let $Z^{0}(t)\equiv Z(t)-tZ(1)$, then:

\begin{enumerate}
\item  Under \textrm{H0 }, we have
\begin{equation*}
n^{-1/2}T\Longrightarrow \sigma _{0}\left( \max_{t\in \left[ 0,1\right]
}W^{0}(t)-\min_{t\in \left[ 0,1\right] }W^{0}(t)\right) ,
\end{equation*}
where $\ \sigma _{0}^{2}=\left( \sum_{k=m}^{\infty }{\frac{{c_{k}^{2}}}{{k!}}%
}\right) ^{-1}\left( \sum_{k=m}^{\infty }{\frac{{c_{k}^{2}}}{{k!}}}%
\sum_{i\in \mathbb{Z}}r^{k}(|i|)\right) $.

\item  Under \textrm{H1, }we have
\begin{equation*}
\frac{1}{d_{n,m}}T\Longrightarrow \sigma _{1}\left( \max_{t\in \left[ 0,1%
\right] }Y_{m}^{0}(t)-\min_{t\in \left[ 0,1\right] }Y_{m}^{0}(t)\right) ,
\end{equation*}
where $\sigma _{1}=\frac{|{c_{m}|}}{m{!}\sqrt{Var(Y)}}.$
\end{enumerate}
\end{corollary}

The proof of Corollary \ref{rss} is easy. It follows directly from the
result of Theorem \ref{mth} and the continuous mapping theorem. We only need
to prove that $S^2$ is a consistent estimator of $\sigma ^{2}=\mathbb{E}%
\left( Y^{2}\right) -\mathbb{E}^{2}\left( Y\right) $. Under H0 and
H1, $r(n)\rightarrow 0$, hence the stationary Gaussian sequence $(X_{i})$ is
mixing, thus ergodic, and so is $(Y_{i})=(G(X_{i}))$. Since
$\mathbb{E}(Y^{2})<\infty $, Birkhoff's ergodic theorem yields
$S^{2}\rightarrow \sigma ^{2}$ almost surely. One of the drawback of this limiting statistics is that two different
constants, $\sigma _{0}$ and $\sigma _{1},$ appeared in the limiting
distributions. Clearly, one can obtain similar results for the modified $R/S$
statistics if we replace $\sigma _{0}$ by some consistent estimator. This is
stated in the following theorem.

\begin{theorem}
\label{th2}  Let $Y=G(X)$, where $(X_{n})$ is a stationary
Gaussian sequence. Let $T(n,q)$ the modified $R/S$ statistics given by :
\begin{equation}
T(n,q)\equiv \frac{R}{\sqrt{\frac{n}{q+1}V_{n}(q)}}  \label{tnq}
\end{equation}
Let $q=q(n)\rightarrow \infty $ and assume that $q=o(n)$
and
\begin{equation}
V_{n}(q)/V(q)\rightarrow 1\ \,\,\quad \text{in probability,}  \label{hos}
\end{equation}
then:

\begin{enumerate}
\item  Under \textrm{H0 }
\begin{equation*}
T(n,q)\Longrightarrow \left(
\max_{t\in \left[ 0,1\right] }W^{0}(t)-\min_{t\in \left[ 0,1\right]
}W^{0}(t)\right) .
\end{equation*}

\item  Under \textrm{H1 }
\begin{equation*}
\left( \frac{n}{q+1}\right) ^{-\frac{1-m\alpha }{2}}\left( \frac{L(n)}{L(q)}%
\right) ^{-m/2}T(n,q)\Longrightarrow \left(
\max_{t\in \left[ 0,1\right] }Y_{m}^{0}(t)-\min_{t\in \left[ 0,1\right]
}Y_{m}^{0}(t)\right) .
\end{equation*}
\end{enumerate}
\end{theorem}

\noindent Under H0, conditions ensuring (\ref{hos}), and even the almost sure
convergence of the Bartlett estimator $S_{n}^{2}(q)$, are given in Theorem A1
of \cite{ber06}. Under H1, the Bartlett (HAC) estimator under long memory is
studied by Abadir, Distaso and Giraitis \cite{adg09} for linear processes; see
also Hosking \cite{hosk} (Theorem 4 p. 267 therein) for the sample
autocovariances of Gaussian long-memory sequences.
Under H1, establishing (\ref{hos}) for non-linear functionals with Hermite
rank $m\geq 2$ requires a law of large numbers for squared block sums, which
is beyond the scope of this paper. Note however that, for testing purposes,
only the consistency of the test under H1 is needed, and this property holds
without assuming (\ref{hos}), see Proposition \ref{prcons} below.

\textbf{Proof.} Under H0, using Theorem \ref{mth} and Slutsky Lemma, we only
need to show that $\frac{V(q)}{q}$ tends to $\sum_{k=m}^{\infty }{\frac{{%
c_{k}^{2}}}{{k!}}}\sum_{i\in \mathbb{Z}}r^{k}(|i|).$ This can be achieved by
the dominated convergence theorem.

Under H1, from Theorem \ref{mth}, we have
\begin{equation*}
\frac{R}{d_{n,m}\bar{\sigma}_{1}}\Longrightarrow \left( \max_{t\in \left[ 0,1%
\right] }Y_{m}^{0}(t)-\min_{t\in \left[ 0,1\right] }Y_{m}^{0}(t)\right)
\end{equation*}
where $\bar{\sigma}_{1}^{2}=(\frac{c_{m}}{m!})^{2}.$ Hence it suffices to
show that
\begin{equation*}
\left( \frac{n}{q+1}\right) ^{-\frac{1-m\alpha }{2}}\left( \frac{L(n)}{L(q)}%
\right) ^{-m/2}\frac{d_{n,m}\bar{\sigma}_{1}}{\sqrt{\frac{n}{q+1}V_{n}(q)}}%
\rightarrow 1,\text{ in probability.}
\end{equation*}
Since $V_{n}(q)/V(q)\rightarrow 1$ in probability and under H1, $V(q)\sim (%
\frac{c_{m}}{m!})^{2}d_{q,m}^{2},$ then
\begin{eqnarray*}
Q(n) &\equiv &\left( \frac{n}{q+1}\right) ^{-\frac{1-m\alpha }{2}}\left(
\frac{L(n)}{L(q)}\right) ^{-m/2}\frac{d_{n,m}\bar{\sigma}_{1}}{\sqrt{\frac{n%
}{q+1}V(q)}} \\
&\sim &\left( \frac{n}{q+1}\right) ^{-\frac{1-m\alpha }{2}}\left( \frac{L(n)%
}{L(q)}\right) ^{-m/2}\frac{d_{n,m}}{\sqrt{\frac{n}{q+1}}\,d_{q,m}}\sim 1.
\end{eqnarray*}
Here $\sim $ denotes equivalence when $n\rightarrow +\infty .$ Again by
Theorem \ref{mth} and Slutsky Lemma, we complete the proof of the Theorem%
\endproof

One of the advantages of these modifications is to overcome the problem of
the different constants $\sigma _{0}$ and $\sigma _{1}$. However, an
important question rises: what is the optimal $q$ to choose? Indeed,
choosing a large $q$ reduces the power of the statistics while choosing a
small $q$ increases its power but decreases the level of the statistics. For
example, taking $q=0$ which corresponds to the original $R/S$ statistics
will typically make the value of statistics high resulting a rejection of H0
even when H0 is true. This fact is well-confirmed by our simulations.

On the other hand, taking a large value of $q$ will typically reduce the
power of the modified statistics resulting a bias towards H0. The problem is
more predominant when the dependence parameter $m\alpha $ is close to 1.
Clearly, according to the second part of Theorem \ref{th2}, when $\frac{n}{%
q+1}$ tends to infinity, the power of the statistics approaches $1$ at the
speed of $\left( \frac{n}{q+1}\right) ^{\frac{1-m\alpha }{2}}.$ So when $%
1-m\alpha $ or $\frac{n}{q+1}$ is small, the value of the statistics will be
typically small resulting in a preference of H0 even when H1 is true. The
presence of this problem is quite natural since when $1-m\alpha $ is small
we are at the phase boundary between long-memory and short-memory. This
makes the test very sensitive and difficult. An alternative to overcome the
problem \ of different constants is, as suggested by Taqqu et al. \cite
{taq99}, to make a regression analysis on $R$ and the sample size in a
log-log scale. However, up to our knowledge, this method is suitable for
estimation of the Hurst parameter and when used for testing it is still lack
of a theoretic study although some empirical studies are available.

\subsection{Consistency under H1 and data-driven choice of $q$}
\label{sec:qdata}

The Bartlett estimator can be written as a mean of squared block sums.
Indeed, let $Z_{i}=Y_{i}-\bar{Y}_{n}$ for $1\leq i\leq n$ and, for
$1-q\leq t\leq n$, let $B_{t}(q)=\sum_{i\in \lbrack t,t+q]\cap \lbrack
1,n]}Z_{i}$. A direct computation gives
\begin{equation}
V_{n}(q)=\frac{1}{n}\sum_{t=1-q}^{n}B_{t}^{2}(q),  \label{vnblock}
\end{equation}
since each pair $(i,j)$ with $|i-j|\leq q$ belongs to exactly $q+1-|i-j|$
blocks. In particular $V_{n}(q)\geq 0$. In the sequel we use Potter's bounds
for slowly varying functions (see \cite{bgt87}, Theorem 1.5.6): for any
$\delta >0$ there exists $C_{\delta }$ such that $L(x)/L(y)\leq C_{\delta
}(x/y)^{\delta }$ for $x\geq y$ large enough.

\begin{lemma}
\label{lemsup} Assume H1 and let $Q=Q_{n}\rightarrow \infty $ with
$Q=o(n)$. Then there exists a constant $K$ such that, for $n$ large enough,
\begin{equation*}
\mathbb{E}\left( \sup_{0\leq q\leq Q}V_{n}(q)\right) \leq K\,V(Q).
\end{equation*}
\end{lemma}

\textbf{Proof.} Let $P_{0}=0$ and $P_{j}=\sum_{i=1}^{j}(Y_{i}-\mathbb{E}(Y))$.
For $q\leq Q$, the set $[t,t+q]\cap \lbrack 1,n]$ is of the form
$[s+1,s+l]$ with $s=\max (t,1)-1$ and $l\leq Q+1$. Hence
\begin{equation*}
|B_{t}(q)|\leq M_{s}+(Q+1)|\bar{Y}_{n}-\mathbb{E}(Y)|,\qquad M_{s}\equiv
\max_{1\leq l\leq Q+1,\,s+l\leq n}|P_{s+l}-P_{s}|.
\end{equation*}
The index $s=0$ corresponds to $Q+1$ values of $t$ and each $s\geq 1$ to one
value of $t$. Let $m(Q)\equiv \mathbb{E}\left( \max_{1\leq l\leq
Q+1}P_{l}^{2}\right) $. By (\ref{vnblock}) and stationarity, we get for
$Q\leq n$
\begin{equation*}
\mathbb{E}\left( \sup_{q\leq Q}V_{n}(q)\right) \leq \frac{2(n+Q)}{n}\left(
m(Q)+(Q+1)^{2}\Var(\bar{Y}_{n})\right) .
\end{equation*}
Under H1, $V(Q)\sim (c_{m}/m!)^{2}d_{Q,m}^{2}$ and
$\Var(\bar{Y}_{n})\leq Kn^{-m\alpha }L^{m}(n)$. Hence, by (\ref{dnm}) and
Potter's bounds,
\begin{equation*}
\frac{(Q+1)^{2}\Var(\bar{Y}_{n})}{V(Q)}\leq K\left( \frac{Q}{n}\right)
^{m\alpha }\left( \frac{L(n)}{L(Q)}\right) ^{m}\leq K\left( \frac{Q}{n}
\right) ^{m(\alpha -\delta )}\rightarrow 0
\end{equation*}
for $\delta <\alpha $. It remains to show that $m(Q)\leq KV(Q)$. Choose
$m^{\ast }$ such that $m^{\ast }\alpha >1$ and write
$Y-\mathbb{E}(Y)=\sum_{k=m}^{m^{\ast }-1}\frac{c_{k}}{k!}H_{k}(X)+G^{m^{\ast
}}(X)$. By Minkowski's inequality and Lemma \ref{lem2} with $p=4$,
\begin{align*}
\left\| \max_{l\leq Q+1}\left| S_{l}(H_{k})\right| \right\| _{2}
&\leq K\left( Q\sum_{|i|\leq Q+1}|r(i)|^{k}\right) ^{1/2} \\
&=O\left( Q^{1-m\alpha /2}L^{m/2}(Q)\right) ,\qquad m\leq k<m^{\ast }.
\end{align*}
For the remainder, since $|\mathrm{Cov}(G^{m^{\ast }}(X_{0}),G^{m^{\ast
}}(X_{i}))|\leq \left\| G\right\| _{2}^{2}|r(i)|^{m^{\ast }}$, the variance
of any sum of $l$ consecutive terms $G^{m^{\ast }}(X_{i})$ is bounded by
$Kl$. By the maximal inequality of Moricz, Serfling and Stout \cite{mor82},
\begin{equation*}
\mathbb{E}\left( \max_{l\leq Q+1}S_{l}^{2}(G^{m^{\ast }})\right) \leq KQ\log
^{2}(2Q)=o(V(Q)),
\end{equation*}
since $V(Q)$ grows like $Q^{2-m\alpha }L^{m}(Q)$ and $2-m\alpha >1$. Hence
$m(Q)\leq KV(Q)$.\endproof

\begin{proposition}
\label{prcons} Assume H1 and let $q=q_{n}\rightarrow \infty $ with
$q=o(n)$. Then $V_{n}(q)=O_{P}(V(q))$ and, for every $\eta >0$, there
exists $c>0$ such that
\begin{equation*}
\liminf_{n\rightarrow \infty }\mathbb{P}\left( \left( \frac{n}{q+1}\right)
^{-\frac{1-m\alpha }{2}}\left( \frac{L(n)}{L(q)}\right) ^{-m/2}T(n,q)\geq
c\right) \geq 1-\eta .
\end{equation*}
In particular $T(n,q)\rightarrow \infty $ in probability, i.e. the test is
consistent.
\end{proposition}

\textbf{Proof.} By Lemma \ref{lemsup} with $Q=q$ and Markov's inequality,
$V_{n}(q)=O_{P}(V(q))$. By Theorem \ref{mth} and the continuous mapping
theorem, $R/d_{n,m}\Longrightarrow \frac{|c_{m}|}{m!}\mathcal{R}$, where
$\mathcal{R}=\max_{t}Y_{m}^{0}(t)-\min_{t}Y_{m}^{0}(t)$. We have
$\mathcal{R}>0$ almost surely. Indeed, $\mathcal{R}=0$ implies
$Y_{m}(1/2)=Y_{m}(1)/2$, while $Y_{m}(1/2)-Y_{m}(1)/2$ is a non-zero element
of the $m$-th Wiener chaos (by self-similarity, its variance equals
$2^{-2H}-1/4>0$, with $H=1-m\alpha /2<1$), whose law is absolutely
continuous by \cite{shi80}. Finally, as in the proof of Theorem \ref{th2},
\begin{equation*}
\frac{d_{n,m}}{\sqrt{\frac{n}{q+1}V(q)}}\sim \frac{m!}{|c_{m}|}\left(
\frac{n}{q+1}\right) ^{\frac{1-m\alpha }{2}}\left( \frac{L(n)}{L(q)}\right)
^{m/2},
\end{equation*}
which tends to infinity by Potter's bounds. \endproof

\medskip

\noindent \textbf{Data-driven choice of }$q$\textbf{.} Proposition
\ref{prcons} shows that the power of the test increases with $n/q$, so that
$q$ should be small when the dependence is strong. Lo \cite{lo91} used the
data-dependent rule (\ref{qlo}) below, based on the rule of Andrews
\cite{and91} which is designed for short memory. As pointed out by Teverovsky
et al. \cite{taq99} and Xiao \cite{xia03}, such data-dependent bandwidths
lead to a strong bias towards H0 under long-memory alternatives; see also
Lima and Xiao \cite{ll10} for a bandwidth which is robust to the presence of
long memory. We propose the following choice:
\begin{equation}
\hat{q}=\left[ n^{\min (1/2,\,1/2-\hat{d})}\right] =\left[
n^{1/2-(\hat{d})_{+}}\right] ,  \label{qhat}
\end{equation}
where $\hat{d}$ is an estimator of the memory parameter $d$ of the sequence
$(Y_{i})$; under H0, $d=0$, and under H1, $d=(1-m\alpha )/2$. For instance,
$\hat{d}$ can be the local Whittle estimator \cite{rob95}, whose consistency
for non-linear transformations of Gaussian processes is established in
\cite{dgh06}, or the Whittle estimator of a parametric ARFIMA$(1,d,0)$ model
\cite{FT86,gt99}, which separates the short-range dependence from the long
memory. The simulations of Section \ref{sec:mc} show that the latter choice
should be preferred.

\begin{proposition}
\label{prqdata} Let $\hat{q}$ be defined by (\ref{qhat}).

\begin{enumerate}
\item  Assume H0 with $\sigma (G)>0$, that $(\hat{d})_{+}\log
n\rightarrow 0$ in probability, and that (\ref{hos}) holds uniformly around
$q_{n}=[n^{1/2}]$, i.e.
$\sup_{q_{n}/2\leq q\leq 2q_{n}}|V_{n}(q)/V(q)-1|\rightarrow 0$ in
probability. Then
\begin{equation*}
T(n,\hat{q})\Longrightarrow \left( \max_{t\in \left[ 0,1\right]
}W^{0}(t)-\min_{t\in \left[ 0,1\right] }W^{0}(t)\right) .
\end{equation*}

\item  Assume H1 and that $\hat{d}\rightarrow d$ in probability. Then
$T(n,\hat{q})\rightarrow \infty $ in probability.
\end{enumerate}
\end{proposition}

\textbf{Proof.} 1. Since $n^{-(\hat{d})_{+}}=\exp (-(\hat{d})_{+}\log
n)\rightarrow 1$ in probability, $\hat{q}/q_{n}\rightarrow 1$ in
probability. Hence $\mathbb{P}(q_{n}/2\leq \hat{q}\leq 2q_{n})\rightarrow 1$
and $V_{n}(\hat{q})/V(\hat{q})\rightarrow 1$ in probability. Under H0,
$V(q)/(q+1)\rightarrow \sigma ^{2}(G)$ as $q\rightarrow \infty $ (see the
proof of Theorem \ref{th2}), so that $V_{n}(\hat{q})/(\hat{q}+1)\rightarrow
\sigma ^{2}(G)$ in probability. Since $T(n,\hat{q})=n^{-1/2}R/(V_{n}(\hat{q}
)/(\hat{q}+1))^{1/2}$, the result follows from Theorem \ref{mth} and
Slutsky's lemma.

2. Let $0<\varepsilon <d/4$ and $Q_{n}=[n^{1/2-d+\varepsilon }]$. Since
$\hat{d}\rightarrow d>0$ in probability, with probability tending to one,
$n^{1/2-d-\varepsilon }-1\leq \hat{q}\leq Q_{n}$, and then
\begin{equation*}
T^{2}(n,\hat{q})\geq \frac{R^{2}}{n^{1/2+d+\varepsilon }\sup_{q\leq
Q_{n}}V_{n}(q)}.
\end{equation*}
By Lemma \ref{lemsup}, $\sup_{q\leq Q_{n}}V_{n}(q)=O_{P}(V(Q_{n}))$, and,
since $2-m\alpha =1+2d$,
\begin{equation*}
V(Q_{n})\leq Kn^{(1/2-d+\varepsilon )(1+2d)}L^{m}(Q_{n}).
\end{equation*} On the other hand, as in the proof of Proposition \ref{prcons},
$R^{2}\geq cd_{n,m}^{2}\geq c^{\prime }n^{1+2d}L^{m}(n)$ with probability
at least $1-\eta $. The exponent of $n$ in the resulting lower bound is
\begin{equation*}
(1+2d)-\left( \frac{1}{2}+d+\varepsilon \right) -\left( \frac{1}{2}
-d+\varepsilon \right) (1+2d)=d+2d^{2}-2\varepsilon (1+d)>0,
\end{equation*}
and the ratio of slowly varying functions is handled by Potter's bounds.
Hence $T(n,\hat{q})\rightarrow \infty $ in probability.\endproof

\medskip

\noindent \textbf{Remarks.} 1. The condition $(\hat{d})_{+}\log
n\rightarrow 0$ under H0 holds for the local Whittle estimator computed with
$M$ frequencies as soon as $\sqrt{M}\hat{d}=O_{P}(1)$ and $(\log n)^{2}/M
\rightarrow 0$, see \cite{rob95} for linear processes and \cite{dgh06} for
non-linear ones. For the Whittle estimator of a correctly specified Gaussian
ARFIMA model, $\hat{d}-d=O_{P}(n^{-1/2})$ \cite{FT86} (see \cite{gt99} for
non-Gaussian sequences), so that the condition holds. Under misspecification,
$\hat{d}$ converges to a pseudo-true value, which may be positive when the
correlations decay slowly (see Section \ref{sec:mc}).

2. Under H0, $\hat{q}$ is close to $n^{1/2}$, while under H1 it behaves like
$n^{1/2-d}$: $q$ decreases when the dependence increases, which improves the
power. In finite samples, a positive bias of $\hat{d}$ under H0 reduces
$\hat{q}$ and may inflate the size of the test, especially when the
correlations decay slowly (see Section 4).

\section{Monte Carlo Simulations}\label{sec:mc}

In this section we study the finite sample performance of the test of H0
(short-range dependence) against H1 (long-range dependence) based on the
modified rescaled range statistic $T(n,q)$, for several choices of the
truncation parameter $q$.

\paragraph{Models.} We consider the following sequences $(Y_{i})$.

\begin{itemize}
\item[M1] Gaussian sequences with covariance function
$r(k)=(1+k^{2})^{-\alpha /2}$, $\alpha \in \{0.2,0.3,\ldots ,1.8\}$. The
sequence is LRD (H1) for $\alpha <1$ and SRD (H0) for $\alpha >1$.

\item[M2] $Y=\exp (X^{2}/5)$, where $X$ is Gaussian with covariance function
$r_{X}(k)=(1+k^{2})^{-\alpha /4}$. Since $G(x)=\exp (x^{2}/5)$ has a finite
second moment and Hermite rank $2$, $\Cov(Y_{0},Y_{k})\sim
c_{2}^{2}/2\,r_{X}^{2}(k)\sim Ck^{-\alpha }$. We use the same grid of values
of $\alpha $ as for M1.

\item[M3] Stationary Gaussian ARFIMA$(0,d,0)$ sequences,
$d\in \{0,0.05,\ldots ,0.45\}$. For $d=0$ the sequence is i.i.d., while for
$d>0$, $r(k)\sim Ck^{2d-1}$ and the sequence is LRD.

\item[M4] $Y=X^{2}$, where $X$ is as in M3. Then
$\Cov(Y_{0},Y_{k})=2r_{X}^{2}(k)\sim Ck^{-\alpha _{d}}$ with
$\alpha _{d}=2-4d$, so that $Y$ is SRD for $d<1/4$ and LRD for $d>1/4$.
This model allows to consider continuously SRD and LRD sequences within the
same model.

\item[M5] Gaussian AR$(1)$ sequences, $r(k)=\phi ^{k}$,
$\phi \in \{0.3,0.5,0.7,0.9\}$. These sequences are SRD but may exhibit a
strong short-range dependence, which is a well-known source of spurious
detection of long memory.
\end{itemize}

The boundary cases $\alpha =1$ (M1, M2) and $d=1/4$ (M4) belong neither to H0
nor to H1 (see Section 3); they are displayed in the figures but are excluded
from the computation of the size. All Gaussian sequences are generated exactly
by circulant embedding (Davies and Harte \cite{dh87}). We consider the sample
sizes $n=500$, $1000$ and $5000$, and also $n=10000$ for M3 and M4. Each
experiment is based on 10\,000 independent replications.

\paragraph{Choices of $q$.} We compare the following truncation parameters:
$q=0$, which corresponds to the classical $R/S$ statistic; the fixed values
$q=[n^{a}]$, $a=0.1,\ldots ,0.5$; the value $q_{Lo}$ suggested by Lo
\cite{lo91}, based on the rule of Andrews \cite{and91},
\begin{equation}
q_{Lo}=[k_{n}],\qquad k_{n}=(3n/2)^{1/3}\left( \frac{2|\hat{\rho}|}{1-\hat{\rho}^{2}}\right) ^{2/3},  \label{qlo}
\end{equation}
where $\hat{\rho}$ is the estimated first-order autocorrelation coefficient of
the data; and the data-driven choice $\hat{q}$ defined in (\ref{qhat}). For
the latter we use two estimators of the memory parameter: the local Whittle
estimator \cite{rob95} computed with $M=[n^{b}]$ frequencies, $b=0.5$, $0.65$
and $0.8$ (denoted by $\hat{q}_{b}$), and the Whittle estimator of an
ARFIMA$(1,d,0)$ model \cite{FT86,gt99} (denoted by $\hat{q}_{W}$); see also
\cite{bbes26} for frequency-domain inference in weak FARIMA models.

\paragraph{Decision rule.} We reject H0 when $T(n,q)>1.747$, which is the
$95\%$ quantile of the limiting statistic under H0,
\begin{equation*}
V\equiv \left( \max_{t\in \left[ 0,1\right] }W^{0}(t)-\min_{t\in \left[ 0,1
\right] }W^{0}(t)\right) .
\end{equation*}
Recall that $E(V)=\sqrt{\pi /2}\approx 1.25$, $E(V^{2})=\pi ^{2}/6$ and
$\sigma (V)\approx 0.27$. Some fractiles of $V$ are given in Table
\ref{tabfractiles}. Since the alternative of interest is long memory, which
increases the range of the partial sums, we use a one-sided test, whereas Lo
\cite{lo91} uses a two-sided one.

\begin{table}[!ht]
\caption{Fractiles of the distribution of $V$ }
\label{tabfractiles}
\begin{center}
\begin{tabular}{lcccccccc}
\hline
$x$ &  & 0.1 & 0.3 & 0.5 & 0.7 & 0.8 & 0.9 & 0.95 \\ \hline
$Q(x)$ &  & 0.927 & 1.09 & 1.223 & 1.374 & 1.473 & 1.620 & 1.747 \\ \hline
\end{tabular}
\end{center}
\end{table}

\begin{figure}[!ht]
\centering
\includegraphics[width=\textwidth]{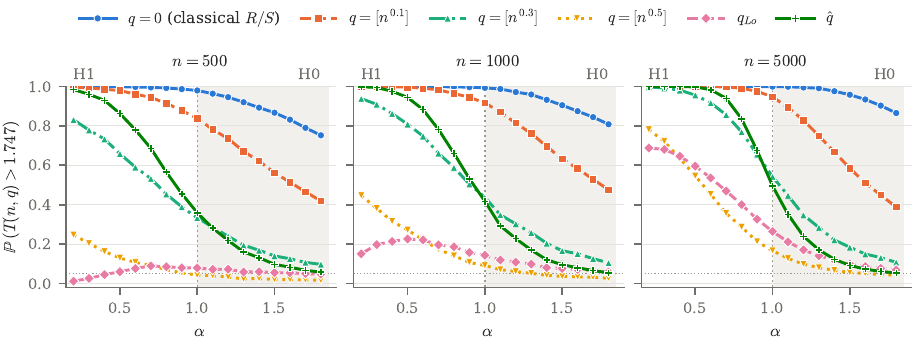}
\caption{Empirical rejection probabilities of the test based on $T(n,q)$
(nominal level 5\%, dotted line) for model M1, as functions of $\alpha$. The
shaded area corresponds to H0 and the dashed line to the boundary
$\alpha =1$.}
\label{figM1}
\end{figure}

\begin{figure}[!ht]
\centering
\includegraphics[width=\textwidth]{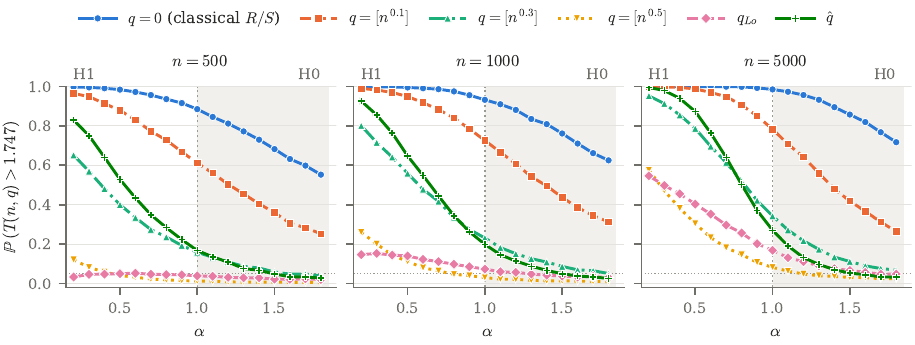}
\caption{Same as Figure \ref{figM1} for model M2, $Y=\exp (X^{2}/5)$.}
\label{figM2}
\end{figure}

\begin{figure}[!ht]
\centering
\includegraphics[width=\textwidth]{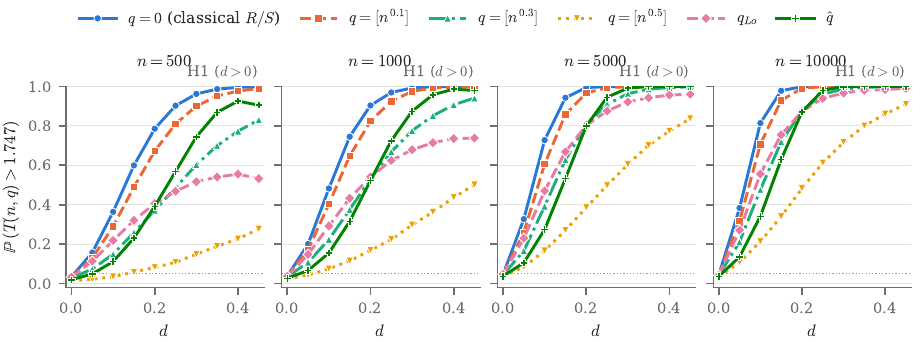}
\caption{Empirical rejection probabilities for model M3 (ARFIMA$(0,d,0)$), as
functions of $d$. The sequence is i.i.d.\ for $d=0$ and LRD for $d>0$.}
\label{figM3}
\end{figure}

\begin{figure}[!ht]
\centering
\includegraphics[width=\textwidth]{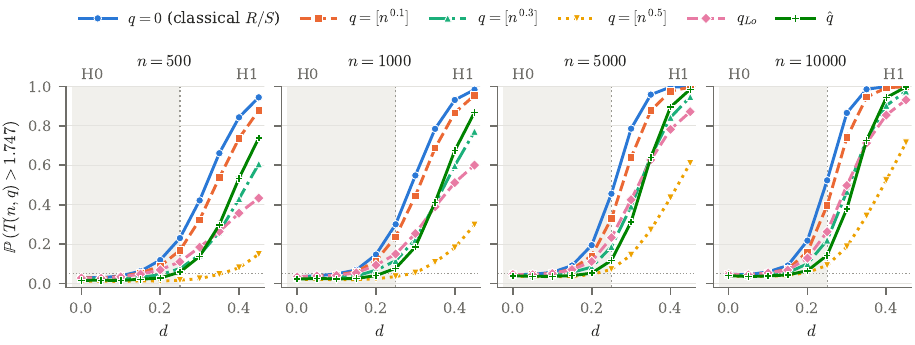}
\caption{Empirical rejection probabilities for model M4 ($Y=X^{2}$, $X$
ARFIMA$(0,d,0)$), as functions of $d$. The shaded area corresponds to H0
($d<1/4$) and the dashed line to the boundary $d=1/4$.}
\label{figM4}
\end{figure}

\begin{table}[!ht]
\centering
\footnotesize
\caption{Empirical size (in \%) under H0. Based on 10000 replications; one-sided test, nominal level 5\% (critical value 1.747). $\hat q_{b}$: data-driven $q$ (\ref{qhat}) with the local Whittle estimator computed on $M=[n^{b}]$ frequencies; $\hat q_{W}$: data-driven $q$ with the Whittle estimator of an ARFIMA$(1,d,0)$ model.}
\label{tabsize}
\setlength{\tabcolsep}{2.4pt}
\begin{tabular}{llrrrrrrrrrrr}
\hline
Model & $n$ & $0$ & $n^{0.1}$ & $n^{0.2}$ & $n^{0.3}$ & $n^{0.4}$ & $n^{0.5}$ & $q_{Lo}$ & $\hat q_{0.5}$ & $\hat q_{0.65}$ & $\hat q_{0.8}$ & $\hat q_{W}$ \\ \hline
M1, $\alpha=1.2$ & 1000 & 97.2 & 81.4 & 53.8 & 30.2 & 14.5 & 6.1 & 10.9 & 24.8 & 30.7 & 63.0 & 22.8 \\
M1, $\alpha=1.2$ & 5000 & 99.5 & 82.9 & 58.6 & 34.8 & 18.6 & 9.9 & 17.1 & 21.3 & 25.3 & 53.9 & 24.0 \\
M1, $\alpha=1.5$ & 1000 & 90.5 & 63.2 & 34.2 & 16.6 & 8.3 & 3.6 & 7.6 & 12.5 & 13.7 & 30.8 & 9.5 \\
M1, $\alpha=1.5$ & 5000 & 95.8 & 58.5 & 34.0 & 18.4 & 9.7 & 5.6 & 9.9 & 10.0 & 10.3 & 20.8 & 9.2 \\
M1, $\alpha=1.8$ & 1000 & 80.8 & 47.6 & 22.8 & 10.5 & 5.6 & 2.6 & 5.8 & 7.9 & 7.6 & 14.7 & 5.4 \\
M1, $\alpha=1.8$ & 5000 & 86.5 & 38.8 & 20.3 & 10.7 & 6.3 & 4.2 & 6.9 & 6.5 & 5.9 & 9.3 & 5.3 \\
\hline
M2, $\alpha=1.2$ & 1000 & 88.1 & 60.8 & 33.1 & 15.1 & 6.1 & 2.1 & 5.7 & 13.4 & 16.1 & 34.0 & 11.4 \\
M2, $\alpha=1.2$ & 5000 & 95.7 & 64.3 & 39.5 & 21.7 & 10.7 & 4.8 & 11.1 & 13.2 & 15.4 & 30.3 & 13.1 \\
M2, $\alpha=1.5$ & 1000 & 76.2 & 43.7 & 20.5 & 8.7 & 3.4 & 1.4 & 3.8 & 7.2 & 7.8 & 16.2 & 5.0 \\
M2, $\alpha=1.5$ & 5000 & 85.8 & 41.9 & 22.2 & 11.2 & 5.7 & 3.1 & 6.6 & 6.3 & 6.3 & 11.5 & 5.3 \\
M2, $\alpha=1.8$ & 1000 & 62.5 & 31.3 & 13.1 & 5.1 & 2.4 & 1.0 & 3.0 & 4.2 & 4.0 & 7.4 & 2.6 \\
M2, $\alpha=1.8$ & 5000 & 71.7 & 26.4 & 13.1 & 6.8 & 4.0 & 2.4 & 4.8 & 4.1 & 3.8 & 5.5 & 3.2 \\
\hline
M3, $d=0$ & 1000 & 3.7 & 3.8 & 3.6 & 3.5 & 3.1 & 2.4 & 3.6 & 3.6 & 3.0 & 2.6 & 2.7 \\
M3, $d=0$ & 10000 & 4.3 & 4.2 & 4.2 & 4.1 & 4.0 & 3.6 & 4.2 & 4.2 & 3.7 & 3.6 & 3.6 \\
\hline
M4, $d=0.1$ & 1000 & 4.9 & 4.3 & 3.9 & 3.4 & 2.9 & 1.9 & 4.1 & 3.4 & 2.6 & 2.2 & 2.3 \\
M4, $d=0.1$ & 10000 & 5.9 & 5.3 & 5.0 & 4.6 & 4.2 & 3.7 & 5.4 & 4.4 & 4.0 & 3.8 & 3.8 \\
M4, $d=0.2$ & 1000 & 14.7 & 11.6 & 8.8 & 6.8 & 4.5 & 2.4 & 9.3 & 5.8 & 4.6 & 3.9 & 4.0 \\
M4, $d=0.2$ & 10000 & 21.8 & 16.1 & 12.8 & 10.0 & 7.5 & 5.6 & 12.9 & 7.4 & 6.7 & 6.5 & 6.5 \\
\hline
M5, $\phi=0.3$ & 1000 & 33.7 & 16.3 & 7.5 & 4.7 & 3.3 & 2.3 & 4.4 & 4.1 & 3.2 & 3.7 & 3.0 \\
M5, $\phi=0.3$ & 5000 & 39.1 & 12.9 & 7.6 & 5.5 & 4.5 & 3.7 & 5.4 & 4.7 & 4.1 & 4.4 & 3.9 \\
M5, $\phi=0.5$ & 1000 & 71.1 & 38.1 & 15.5 & 6.9 & 3.9 & 2.3 & 4.2 & 5.1 & 4.5 & 10.3 & 3.7 \\
M5, $\phi=0.5$ & 5000 & 77.1 & 27.6 & 12.5 & 7.5 & 5.3 & 4.1 & 5.7 & 5.3 & 4.8 & 6.8 & 4.5 \\
M5, $\phi=0.7$ & 1000 & 97.1 & 75.6 & 37.6 & 13.7 & 5.4 & 2.4 & 3.5 & 8.4 & 15.4 & 75.1 & 6.6 \\
M5, $\phi=0.7$ & 5000 & 98.4 & 60.9 & 27.6 & 10.4 & 5.0 & 3.5 & 4.5 & 5.0 & 5.3 & 53.8 & 4.2 \\
M5, $\phi=0.9$ & 1000 & 100.0 & 99.8 & 91.2 & 57.8 & 20.8 & 4.8 & 1.3 & 65.1 & 99.7 & 99.8 & 9.3 \\
M5, $\phi=0.9$ & 5000 & 100.0 & 99.3 & 85.5 & 43.5 & 13.1 & 5.0 & 4.1 & 17.2 & 96.2 & 100.0 & 5.8 \\
\hline
\end{tabular}
\end{table}

\begin{table}[!ht]
\centering
\footnotesize
\caption{Empirical power (in \%) under H1. Based on 10000 replications; one-sided test, nominal level 5\% (critical value 1.747). $\hat q_{b}$: data-driven $q$ (\ref{qhat}) with the local Whittle estimator computed on $M=[n^{b}]$ frequencies; $\hat q_{W}$: data-driven $q$ with the Whittle estimator of an ARFIMA$(1,d,0)$ model.}
\label{tabpower}
\setlength{\tabcolsep}{2.4pt}
\begin{tabular}{llrrrrrrrrrrr}
\hline
Model & $n$ & $0$ & $n^{0.1}$ & $n^{0.2}$ & $n^{0.3}$ & $n^{0.4}$ & $n^{0.5}$ & $q_{Lo}$ & $\hat q_{0.5}$ & $\hat q_{0.65}$ & $\hat q_{0.8}$ & $\hat q_{W}$ \\ \hline
M1, $\alpha=0.5$ & 1000 & 100.0 & 99.7 & 95.3 & 80.8 & 56.2 & 27.2 & 22.6 & 85.2 & 96.0 & 99.7 & 94.5 \\
M1, $\alpha=0.5$ & 5000 & 100.0 & 100.0 & 99.6 & 95.7 & 80.3 & 54.6 & 59.6 & 95.9 & 99.6 & 100.0 & 99.8 \\
M1, $\alpha=0.8$ & 1000 & 99.8 & 96.5 & 82.5 & 58.9 & 34.6 & 14.1 & 18.6 & 58.1 & 74.6 & 96.0 & 66.0 \\
M1, $\alpha=0.8$ & 5000 & 100.0 & 99.1 & 93.4 & 77.0 & 51.5 & 28.7 & 40.0 & 67.8 & 82.3 & 98.7 & 83.6 \\
\hline
M2, $\alpha=0.5$ & 1000 & 99.5 & 94.9 & 79.5 & 55.9 & 30.6 & 11.2 & 14.0 & 56.2 & 71.9 & 92.7 & 64.8 \\
M2, $\alpha=0.5$ & 5000 & 100.0 & 98.8 & 93.4 & 78.7 & 55.2 & 30.4 & 40.3 & 72.9 & 85.2 & 98.2 & 87.3 \\
M2, $\alpha=0.8$ & 1000 & 97.5 & 84.1 & 58.8 & 33.6 & 15.5 & 4.7 & 9.8 & 31.9 & 42.6 & 71.4 & 34.4 \\
M2, $\alpha=0.8$ & 5000 & 99.8 & 90.6 & 74.4 & 52.0 & 30.3 & 14.6 & 25.7 & 40.7 & 51.2 & 78.7 & 50.2 \\
\hline
M3, $d=0.1$ & 1000 & 48.3 & 40.3 & 31.2 & 22.4 & 14.4 & 7.4 & 29.2 & 19.4 & 16.7 & 14.8 & 15.5 \\
M3, $d=0.1$ & 10000 & 81.4 & 70.7 & 60.0 & 47.7 & 34.1 & 21.3 & 55.6 & 36.0 & 34.8 & 34.0 & 34.1 \\
M3, $d=0.3$ & 1000 & 99.3 & 97.5 & 91.1 & 77.5 & 56.2 & 29.7 & 67.9 & 80.7 & 86.1 & 87.4 & 87.5 \\
M3, $d=0.3$ & 10000 & 100.0 & 100.0 & 99.9 & 98.7 & 92.0 & 71.7 & 96.7 & 99.0 & 99.8 & 99.9 & 99.9 \\
\hline
M4, $d=0.35$ & 1000 & 78.5 & 68.7 & 55.6 & 40.5 & 25.3 & 10.8 & 39.0 & 39.2 & 41.1 & 41.0 & 41.2 \\
M4, $d=0.35$ & 10000 & 98.6 & 94.6 & 86.2 & 73.0 & 54.8 & 34.3 & 71.3 & 65.0 & 69.7 & 72.7 & 72.7 \\
M4, $d=0.45$ & 1000 & 98.4 & 95.5 & 89.1 & 77.2 & 57.1 & 29.8 & 60.1 & 79.3 & 85.2 & 88.3 & 86.9 \\
M4, $d=0.45$ & 10000 & 100.0 & 100.0 & 99.7 & 97.6 & 90.0 & 71.5 & 93.2 & 97.2 & 99.3 & 99.8 & 99.8 \\
\hline
\end{tabular}
\end{table}

\begin{table}[!ht]
\centering
\footnotesize
\caption{Mean and standard deviation of $T(n,q)$, model M1, $n=1000$. Under H0 the limit $V$ has mean 1.25 and standard deviation 0.27. Based on 10000 replications.}
\label{tabmoments}
\setlength{\tabcolsep}{2.4pt}
\begin{tabular}{lrrrrrrrrrrr}
\hline
$\alpha$ & $0$ & $n^{0.1}$ & $n^{0.2}$ & $n^{0.3}$ & $n^{0.4}$ & $n^{0.5}$ & $q_{Lo}$ & $\hat q_{0.5}$ & $\hat q_{0.65}$ & $\hat q_{0.8}$ & $\hat q_{W}$ \\ \hline
0.5 \quad mean & 5.06 & 3.75 & 2.86 & 2.25 & 1.84 & 1.54 & 1.51 & 2.99 & 3.37 & 3.75 & 3.10 \\
\phantom{0.0} \quad s.d. & 1.33 & 0.97 & 0.72 & 0.55 & 0.42 & 0.31 & 0.29 & 1.19 & 1.06 & 0.97 & 0.99 \\[2pt]
0.8 \quad mean & 3.88 & 2.93 & 2.30 & 1.89 & 1.61 & 1.42 & 1.46 & 2.04 & 2.24 & 2.91 & 2.05 \\
\phantom{0.0} \quad s.d. & 1.00 & 0.75 & 0.58 & 0.46 & 0.37 & 0.29 & 0.31 & 0.77 & 0.69 & 0.75 & 0.60 \\[2pt]
1.2 \quad mean & 2.89 & 2.25 & 1.84 & 1.58 & 1.41 & 1.31 & 1.37 & 1.51 & 1.58 & 1.97 & 1.49 \\
\phantom{0.0} \quad s.d. & 0.71 & 0.54 & 0.44 & 0.37 & 0.31 & 0.26 & 0.29 & 0.44 & 0.41 & 0.50 & 0.37 \\[2pt]
1.5 \quad mean & 2.45 & 1.95 & 1.63 & 1.44 & 1.32 & 1.26 & 1.31 & 1.36 & 1.39 & 1.59 & 1.34 \\
\phantom{0.0} \quad s.d. & 0.58 & 0.46 & 0.38 & 0.32 & 0.29 & 0.25 & 0.28 & 0.34 & 0.33 & 0.38 & 0.30 \\[2pt]
1.8 \quad mean & 2.19 & 1.76 & 1.51 & 1.36 & 1.28 & 1.24 & 1.29 & 1.31 & 1.31 & 1.41 & 1.28 \\
\phantom{0.0} \quad s.d. & 0.50 & 0.40 & 0.34 & 0.30 & 0.27 & 0.24 & 0.27 & 0.29 & 0.28 & 0.32 & 0.26 \\[2pt]
\hline
\end{tabular}
\end{table}

\paragraph{Results.} Figures \ref{figM1}--\ref{figM4} display the empirical
rejection probabilities as functions of the dependence parameter, Table
\ref{tabsize} reports the empirical size under H0 and Table \ref{tabpower}
the empirical power under H1. The mean and standard deviation of $T(n,q)$
are reported in Table \ref{tabmoments} for model M1. We highlight the
following conclusions.

\begin{enumerate}
\item The classical $R/S$ statistic ($q=0$) has a very strong preference
towards H1. For instance, for model M1 with $\alpha =1.5$ and $n=1000$, H0 is
rejected in 90.5\% of the cases, and in 80.8\%
of the cases for $\alpha =1.8$. It is therefore necessary to include the
cross-correlations in the normalization, as proposed by Lo \cite{lo91}.

\item For fixed values of $q$, the choice of $q$ results in a trade-off
between size and power. Small values of $q$ give a high power but a large
size near the boundary, while $q=[n^{0.5}]$ keeps the size under control but
loses most of the power: for model M1 with $\alpha =0.8$ and $n=5000$, the
power is only 28.7\%.

\item The choice $q_{Lo}$ controls the size, even for strongly dependent
AR$(1)$ sequences, but it has a low power under long-memory alternatives
(40.0\% for model M1 with $\alpha =0.8$ and $n=5000$). This is
in line with the findings of Teverovsky et al. \cite{taq99} and Xiao
\cite{xia03}. Note that $q_{Lo}$ increases with $\hat{\rho}$, that is with the
strength of the dependence, which reduces the power.

\item The data-driven choice $\hat{q}_{b}$, based on the local Whittle
estimator, has a high power, but it is severely oversized for strongly
dependent AR$(1)$ sequences: for $\phi =0.9$, H0 is rejected in
99.7\% of the cases for $n=1000$ and
96.2\% for $n=5000$ (with $b=0.65$). Indeed, the
local Whittle estimator does not separate the short-range dependence from
the long memory, and the resulting positive estimate of $d$ reduces $q$.

\item The data-driven choice $\hat{q}_{W}$, based on the Whittle estimator of
an ARFIMA$(1,d,0)$ model, combines the advantages of both approaches. Its
size is close to, or slightly below, the nominal level for i.i.d.\
sequences (3.6\% for $n=10000$), for model M4 under
H0 and for AR$(1)$ sequences (5.8\% for $\phi =0.9$
and $n=5000$), although some size distortion remains for strongly dependent
AR$(1)$ sequences in small samples (9.3\% for
$\phi =0.9$ and $n=1000$). Its power is comparable to the one of $\hat{q}_{b}$:
83.6\% for model M1 with $\alpha =0.8$ and $n=5000$,
and 87.5\% (against 67.9\% for
$q_{Lo}$) for model M3 with $d=0.3$ and $n=1000$.

\item The most difficult situation is the vicinity of the boundary between
short and long memory. For models M1 and M2 with $\alpha $ slightly larger
than $1$, all the small values of $q$ are oversized: for model M1 with
$\alpha =1.2$ and $n=5000$, the size is 24.0\% for
$\hat{q}_{W}$, 17.1\% for $q_{Lo}$ and
9.9\% for $q=[n^{0.5}]$. In this case the correlations
decay like $k^{-1.2}$, which is hardly distinguishable from long memory in
finite samples, and the ARFIMA$(1,d,0)$ model is misspecified: the mean of
the estimated memory parameter is 0.13. This illustrates the
phase boundary problem discussed in Section 3.
\end{enumerate}

To summarize, the choice of $q$ is crucial and has to be adapted to the
strength of the dependence: $q$ should be small when the dependence is strong
and close to $n^{1/2}$ under short memory. Among the rules considered here,
$\hat{q}_{W}$ provides the best compromise between size and power, and we
recommend its use in practice.

\section{Proofs.}

In this section we present the proofs of our main results. The main
technical difficulty is to prove the tightness of the Donsker's line under
the minimal conditions. Indeed, if we assume $\mathbb{E}(G^{p}(X))<\infty $
for $p>2$, the proof is straightforward by using the maximal inequality
of \cite[Theorem 2.3]{bb25}. However, in this paper we only assume that the
second moment is finite. Then we need to use a truncation argument to reduce
the moment condition to the minimal one, namely $\mathbb{E}(G^{2}(X)) $ $\ $%
is finite. This makes the proof more technical, since in that case the
Hermite rank of the truncated function $G_{M}\equiv G1_{|G|\leq M}$ may be
smaller \ than the Hermite rank of G. To overcome this difficulty we apply
Hermite expansion to $G_{M}$ and deal with the first few terms separately.
Our proof is based on the following two Lemmas. In the sequel, $K$ denotes
generic constant which value may change in each appearance.

\begin{lemma}[{\cite[Theorem 2.3]{bb25}}]
\label{lem1} Let $m\geq 1$ and assume that $r^{m}\in \mathbb{L}^{1}$. Then
there exists a constant $K=K(r,m)$ such that, for any function $G$ with
Hermite rank greater than or equal to $m$ and
$\mathbb{E}\left[ G^{4}(X)\right] <\infty $, and for all $n\geq 1$,
\begin{equation}
\mathbb{E}\left[ \max_{k\leq n}\left| S_{k}\left( G-\mathbb{E}(G(X))\right)
\right| \right] ^{4}\leq K(r,m)\left( \left( \sqrt{n}\left\| G\left( X\right)
\right\| _{2}\right) ^{4}+n\left\| G(X)\right\| _{4}^{4}\right) .
\label{inef1}
\end{equation}
\end{lemma}

This is Theorem 2.3 in \cite{bb25}, where it is derived from the moment
inequality of Proposition 2.2 therein. The following lemma is useful when the
expansion of $G$ is finite.

\begin{lemma}[{\cite[Theorem 2.1]{bb25}}]
\label{lem2} Let $G$ be a function with Hermite rank $m$ such that
$\mathbb{E}\left[ G\left( X\right) \right] =0$,
$\mathbb{E}\left[ G^{2}\left( X\right) \right] <\infty $, and let $p>2$ be an
integer. Then
\begin{equation*}
\mathbb{E}\left( \sum_{i=1}^{n}G(X_{i})\right) ^{p}\leq \left(
\sum_{k=m}^{\infty }{\frac{{|c_{k}|}}{\sqrt{k!}}}\left(
2n\sum_{i=-n}^{n}|r^{k}(i)|\right) ^{1/2}(p-1)^{k/2}\right) ^{p},
\end{equation*}
and, for every $l\geq 1$, there exists a constant $K=K(p,l)$ independent of
$n$ such that
\begin{equation}
\mathbb{E}\max_{k\leq n}\left| \sum_{i=1}^{k}H_{l}(X_{i})\right| ^{p}\leq
K\left( n\sum_{i=-n}^{n}|r^{l}(i)|\right) ^{p/2}.  \label{lem2eq2}
\end{equation}
\end{lemma}

The first part of Lemma \ref{lem2} is the discrete analogue of Proposition 2
in \cite{ben02}; the second part follows from the first one and Theorem 3.1
in \cite{mor82}, see \cite[Theorem 2.1]{bb25}. Note that Lemma \ref{lem2}
gives sharp bounds for both short-range and long-range dependent processes.
Indeed, the term ${\frac{{c_{k}^{2}}}{{k!}}}\, n\sum_{i=-n}^{n}|r^{k}(i)|$ is
an upper bound for the variance of
$\sum_{i=1}^{n}\frac{c_{k}}{k!}H_{k}(X_{i})$, so the $p$-th moment of the sum
is controlled by the variance to the power $p/2$. Lemmas \ref{lem1} and
\ref{lem2}, together with the Hermite expansion, allow us to control the
fourth moment for both short and long-range dependent sequences.

\subsection{Proof of Theorem \ref{mth}}

We first prove (\ref{res1}), namely the remainder term in our expansion
converges to the Brownian motion. The proof for finite dimensional
convergence is straightforward ( for example see \cite{ben02} or \cite{bre83}%
). Therefore we only have to prove the tightness or the stochastic
equicontinuity.

\subsubsection*{\textbf{Proof of tightness}}

We define
\begin{eqnarray*}
Z_{n}\left( G\right) &\equiv &n^{-1/2}S_{n}\left( G-E(G(X))\right) , \\
w\left( n,\delta \right) &\equiv &\sup_{\left| t-s\right| \leq \delta
}\left| Z_{\left[ nt\right] }\left( G\right) -Z_{\left[ ns\right] }\left(
G\right) \right| .
\end{eqnarray*}
Let $\ M\equiv \sqrt{n\delta }\left\| G\right\| _{2},G_{M}\equiv G1_{\left|
G\right| \leq M},$ $\tilde{G}_{M}\equiv G-G_{M}.$ Then we can write:
\begin{eqnarray*}
\mathbb{P}\left( w\left( n,\delta \right) >\varepsilon \right) &\leq
&\sum_{i=0}^{\left[ 1/\delta \right] }\mathbb{P}\left( \max_{in\delta \leq
k<(i+1)n\delta }\left| Z_{k}\left( G\right) -Z_{\left[ in\delta \right]
}\left( G\right) \right| >\frac{\varepsilon }{3}\right) \\
&\leq &\sum_{i=0}^{\left[ 1/\delta \right] }\mathbb{P}\left( \max_{in\delta
\leq k<(i+1)n\delta }\left| Z_{k}\left( G_{M}\right) -Z_{\left[ in\delta %
\right] }\left( G_{M}\right) \right| >\frac{\varepsilon }{6}\right) \\
&&+\sum_{i=0}^{\left[ 1/\delta \right] }\mathbb{P}\left( \max_{in\delta \leq
k<(i+1)n\delta }\left| Z_{k}\left( \tilde{G}_{M}\right) -Z_{\left[ in\delta %
\right] }\left( \tilde{G}_{M}\right) \right| >\frac{\varepsilon }{6}\right)
\\
&\equiv &E1+E2.
\end{eqnarray*}

\textbf{Control of E1:} Since the underlying sequence is stationary, we have

\begin{eqnarray*}
E1 &=&\sum_{i=0}^{\left[ 1/\delta \right] }\mathbb{P}\left( \max_{in\delta
\leq k<(i+1)n\delta }\left| Z_{k}\left( G_{M}\right) -Z_{\left[ in\delta %
\right] }\left( G_{M}\right) \right| >\frac{\varepsilon }{6}\right) \\
&\leq &\left( \left[ 1/\delta \right] +1\right) \mathbb{P}\left( \max_{0\leq
k<n\delta }\left| Z_{k}\left( G_{M}\right) \right| >\frac{\varepsilon }{6}%
\right) \\
&\leq &\frac{2}{\delta }\mathbb{P}\left( \max_{0\leq k<n\delta }\left|
Z_{k}\left( G_{M}\right) \right| >\frac{\varepsilon }{6}\right) .
\end{eqnarray*}

\bigskip Now we expand $G_{M}$ in terms of Hermite polynomials to obtain
\begin{eqnarray*}
G_{M} &=&\sum\limits_{l=0}^{\infty }{\frac{c_{l\left( M\right) }}{l{!}}}H_{l}
\\
&=&\sum\limits_{l=0}^{m-1}{\frac{c_{l\left( M\right) }}{l{!}}}%
H_{l}+\sum\limits_{l=m}^{\infty }{\frac{c_{l\left( M\right) }}{l{!}}}H_{l} \\
&\equiv &\sum\limits_{l=0}^{m-1}{\frac{c_{l\left( M\right) }}{l{!}}}%
H_{l}+G_{M}^{m}.
\end{eqnarray*}
Hence,
\begin{eqnarray*}
E_{1} &\leq &\frac{2}{\delta }\sum\limits_{l=1}^{m-1}\mathbb{P}\left(
\max_{0\leq k<n\delta }\left| Z_{k}\left( {\frac{c_{l}\left( M\right) }{l{!}}%
}H_{l}\right) \right| >\frac{\varepsilon }{6m}\right) \\
&&+\frac{2}{\delta }\mathbb{P}\left( \max_{0\leq k<n\delta }\left|
Z_{k}\left( G_{M}^{m}\right) \right| >\frac{\varepsilon }{6m}\right) \\
&\leq &\frac{2}{\delta }\sum\limits_{l=1}^{m-1}\left( {\frac{6mc_{l}\left(
M\right) }{\varepsilon l{!}}}\right) ^{4}\mathbb{E}\left( \max_{0\leq
k<n\delta }\left| Z_{k}\left( H_{l}\right) \right| \right) ^{4} \\
&&+\frac{2}{\delta }\left( \frac{6m}{\varepsilon }\right) ^{4}\mathbb{E}%
\left( \max_{0\leq k<n\delta }\left| Z_{k}\left( G_{M}^{m}\right) \right|
\right) ^{4} \\
&\equiv &\frac{2}{\delta }\sum\limits_{l=1}^{m-1}E_{1,l}+\frac{2}{\delta }%
E_{1,m}
\end{eqnarray*}
For $0<l<m$, we apply Lemma \ref{lem2} with p=4 to get
\begin{eqnarray}
\mathbb{E}\left( \max_{0\leq k<n\delta }\left| Z_{k}\left( H_{l}\right)
\right| \right) ^{4} &\leq &\delta ^{2}K\left( \sum_{i=0}^{n\delta
}|r(i)|^{l}\right) ^{2}  \notag \\
&\leq &\delta ^{2}K\left( \sum_{i=0}^{n}|r(i)|^{m}\right) ^{2l/m}\left(
\delta n\right) ^{2(1-l/m)}.  \label{eq33}
\end{eqnarray}
Now we control $c_{l}\left( M\right) .$ Since $c_{l}=0,$ for $l<m$, we have
\begin{equation*}
c_{l}\left( M\right) =\mathbb{E}\left( G_{M}\left( X\right) H_{l}\left(
X\right) \right) =-\mathbb{E}\left( \tilde{G}_{M}\left( X\right) H_{l}\left(
X\right) \right) .
\end{equation*}
Using H\"{o}lder's inequality, for $p$ and $q$ such that $p^{-1}+q^{-1}=1,$ $%
1<p\leq 2,$
\begin{eqnarray}
\left| c_{l}\left( M\right) \right| &\leq &\left\| \tilde{G}_{M}\left(
X\right) \right\| _{p}\left\| H_{l}\left( X\right) \right\| _{q}  \notag \\
&\leq &M^{1-2/p}\left\| \tilde{G}_{M}\left( X\right) \right\|
_{2}^{2/p}\left\| H_{l}\left( X\right) \right\| _{q}  \notag \\
&\leq &M^{1-2/p}\left\| G\left( X\right) \right\| _{2}^{2/p}\left\|
H_{l}\left( X\right) \right\| _{q}.  \label{eq34}
\end{eqnarray}
Combining (\ref{eq33}) and (\ref{eq34}), we get for $l<m$
\begin{eqnarray*}
E_{1,l} &\leq &K(\varepsilon ,l,m)\left( M^{1-2/p}\left\| G\left( X\right)
\right\| _{2}^{2/p}\left\| H_{l}\left( X\right) \right\| _{q}\right) ^{4} \\
&&\times \,\delta ^{2}K\left( \sum_{i=0}^{n}|r(i)|^{m}\right) ^{2l/m}\left(
\delta n\right) ^{2(1-l/m)} \\
&\leq &K(\varepsilon ,m,r,G)n^{-2\frac{-2pm+2m+pl}{pm}}\delta ^{-2\frac{%
-3pm+2m+pl}{pm}}.
\end{eqnarray*}
Choosing $p$ such that $-2pm+2m+pl>0,$ i.e. $p<\frac{2m}{2m-l}$, gives
\begin{equation}
\lim_{n\rightarrow \infty }E_{1,l}=0.  \label{be1l}
\end{equation}

Finally, in order to control $E_{1,m}$, we apply Lemma \ref{lem1} to obtain
\begin{equation*}
\mathbb{E}\left( \max_{0\leq k<n\delta }\left| Z_{k}\left( G_{M}^{m}\right)
\right| \right) ^{4}\leq Kn^{-2}\left[ \left( \sqrt{n\delta }\left\|
G_{M}^{m}\right\| _{2}\right) ^{4}+n\delta \left\| G_{M}^{m}\right\| _{4}^{4}%
\right] .
\end{equation*}
Since $\mathbb{E}\left[ G_{M}^{4}(X)\right] \leq M^{2}\,\mathbb{E}\left[
G^{2}(X)\right] $ and $|c_{l}(M)|\leq \left\| G\right\| _{2}\left\|
H_{l}\right\| _{2}$, we have
\begin{equation*}
\left\| G_{M}^{m}\right\| _{4}^{4}\leq K\left( M^{2}\left\| G\right\|
_{2}^{2}+\left\| G\right\| _{2}^{4}\right) ,
\end{equation*}
so that, with $M=\sqrt{n\delta }\left\| G\right\| _{2}$ and $n\delta \geq 1$,
$n\delta \left\| G_{M}^{m}\right\| _{4}^{4}\leq K\left( \sqrt{n\delta
}\left\| G\right\| _{2}\right) ^{4}$. Hence
\begin{equation*}
\mathbb{E}\left( \max_{0\leq k<n\delta }\left| Z_{k}\left( G_{M}^{m}\right)
\right| \right) ^{4}\leq Kn^{-2}\left( \sqrt{n\delta }\left\| G\right\|
_{2}\right) ^{4}.
\end{equation*}
Note that $M$ is chosen in such a way that the two terms in the right hand
side of the moment inequality (\ref{inef1}) are of the same order of
magnitude, namely $M=\sqrt{n\delta }\left\| G\right\| _{2},$ hence
\begin{equation}
\mathbb{E}\left( \max_{0\leq k<n\delta }\left| Z_{k}\left( G_{M}^{m}\right)
\right| \right) ^{4}\leq K\left( \left( \sqrt{\delta }\left\| G\right\|
_{2}\right) ^{4}\right) .  \label{bE1m}
\end{equation}
From (\ref{be1l}) (\ref{bE1m}) we get
\begin{equation*}
\lim_{\delta \rightarrow 0}\lim_{n\rightarrow \infty }E_{1}=0.
\end{equation*}

\textbf{Control of E2: }Since $\left\| G\right\| _{2}$ is finite, we
conclude that
\begin{equation*}
\lim_{n\rightarrow \infty }\mathbb{E}\left( G^{2}1_{|G|\geq \sqrt{n\delta }%
\left\| G\right\| _{2}}\right) =0,
\end{equation*}
for every $\delta >0.$ From \ the inequality
\begin{equation*}
\mathbb{E}\left( G^{2}1_{|G|\geq \sqrt{n\delta }\left\| G\right\|
_{2}}\right) \geq \sqrt{n\delta }\left\| G\right\| _{2}\mathbb{E}\left(
\left| G\right| 1_{|G|\geq \sqrt{n\delta }\left\| G\right\| _{2}}\right) ,
\end{equation*}
we \ deduce that, for any $\delta >0,$%
\begin{equation*}
\sqrt{n}\mathbb{E}\left( \left| \tilde{G}_{M}\right| \right) \leq \left(
\sqrt{\delta }\left\| G\right\| _{2}\right) ^{-1}\mathbb{E}\left(
G^{2}1_{|G|\geq \sqrt{n\delta }\left\| G\right\| _{2}}\right) \rightarrow 0.
\end{equation*}
By Markov's inequality

\begin{eqnarray*}
E_{2} &\equiv &\sum_{i=0}^{\left[ 1/\delta \right] }\mathbb{P}\left(
\max_{in\delta \leq k<(i+1)n\delta }\left| Z_{k}\left( \tilde{G}_{M}\right)
-Z_{\left[ in\delta \right] }\left( \tilde{G}_{M}\right) \right| >\frac{%
\varepsilon }{6}\right) \\
&\leq &\frac{6}{\varepsilon }\sum_{i=0}^{\left[ 1/\delta \right] }\mathbb{E}%
\left( \max_{in\delta \leq k<in\delta +n\delta }\left| Z_{k}\left( \tilde{G}%
_{M}\right) -Z_{\left[ in\delta \right] }\left( \tilde{G}_{M}\right) \right|
\right) \\
&\leq &\frac{K}{\varepsilon \delta }\sqrt{n}\delta \mathbb{E}\left( \left|
\tilde{G}_{M}\right| \right) .
\end{eqnarray*}
Therefore we have $\lim_{n\rightarrow \infty }E_{2}=0$. This concludes the
proof of tightness and the proof of (\ref{res1}). The proof of (\ref{res2})
is an immediate consequence of (\ref{res1}).\endproof

\subsubsection*{Proof of (\ref{res3}).}

Let $\left\| .\right\| $\ denotes the sup-norm in the space $\mathcal{D}%
\left( \left[ 0,1\right] \right) .$ We observe that
\begin{eqnarray*}
D_{n}(G,H_{m}) &\equiv &\mathbb{E}\left( \left\| S_{n}(\cdot ,G-\mathbb{E}%
\left[ G\left( X\right) \right] )-{\frac{{c_{m}}}{m{!}}}S_{n}(\cdot
,H_{m})\right\| \right) \\
&\leq &\sum_{k=m+1}^{m^{\ast }-1}{\frac{\left| {c_{k}}\right| }{k{!}}}%
\mathbb{E}\left( \left\| S_{n}(\cdot ,H_{k})\right\| \right) +\mathbb{E}%
\left( \left\| S_{n}(\cdot ,G^{m^{\ast }})\right\| \right) \\
&\leq &\sum_{k=m+1}^{m^{\ast }-1}{\frac{\left| {c_{k}}\right| }{k{!}}}\left(
\mathbb{E}\left( \max_{j\leq n}\left| S_{j}(H_{k})\right| \right) ^{4}\right)
^{1/4}+\mathbb{E}\left( \left\| S_{n}(\cdot ,G^{m^{\ast }})\right\| \right) .
\end{eqnarray*}
By Lemma \ref{lem2} with $p=4$, for $m+1\leq k<m^{\ast }$,
\begin{equation*}
\left( \mathbb{E}\left( \max_{j\leq n}\left| S_{j}(H_{k})\right| \right)
^{4}\right) ^{1/4}\leq K\left( n\sum_{i=-n}^{n}|r(i)|^{k}\right) ^{1/2},
\end{equation*}
and
\begin{equation*}
\left( n\sum_{i=-n}^{n}|r(i)|^{k}\right) ^{1/2}\leq K\left\{
\begin{array}{ll}
n^{1-k\alpha /2}L^{k/2}(n) & \text{if }k\alpha <1, \\
\left( n\,\ell _{k}(n)\right) ^{1/2} & \text{if }k\alpha =1, \\
n^{1/2} & \text{if }k\alpha >1,
\end{array}
\right.
\end{equation*}
where $\ell _{k}(n)=\sum_{i=1}^{n}i^{-1}L^{k}(i)$ is slowly varying. Choose
$m^{\ast }$ such that $m^{\ast }\alpha >1$. From Lemma \ref{lem1} and arguing
as in the proof of tightness with $\delta =1$, we obtain
\begin{equation*}
\mathbb{E}\left( \left\| S_{n}\left( \cdot ,G^{m^{\ast }}\right) \right\|
\right) \leq K\left( m^{\ast },G,r\right) n^{1/2}.
\end{equation*}
Therefore, there exists a slowly varying function $\tilde{L}$ such that
\begin{eqnarray}
\lefteqn{\mathbb{E}\left( \left\| S_{n}(\cdot ,G-\mathbb{E}\left[ G\left(
X\right) \right] )-\frac{c_{m}}{m!}S_{n}(\cdot ,H_{m})\right\| \right) }
\notag \\
&\leq &K\left( m^{\ast },G,r\right) \left( n^{1-\left( m+1\right) \alpha
/2}+n^{1/2}\right) \tilde{L}(n).  \label{fin}
\end{eqnarray}
From (\ref{dnm}) and (\ref{fin}) we obtain
\begin{eqnarray*}
\lefteqn{\mathbb{E}\left( \frac{1}{d_{n,m}}\left\| S_{n}(\cdot ,G-\mathbb{E}\left[
G\left( X\right) \right] )-{\frac{{c_{m}}}{m{!}}}S_{n}(\cdot ,H_{m})\right\|
\right) } \\
&\leq &K\frac{\tilde{L}(n)}{L^{m/2}(n)}\left( n^{-\alpha /2}+n^{-\left(
1-m\alpha \right) /2}\right) .
\end{eqnarray*}
Since the right hand side of the above inequality tends to zero as $n$ tends
to $\infty ,$ $\frac{1}{d_{n,m}}S_{n}(\cdot ,G-\mathbb{E}\left[ G\left(
X\right) \right] )$ and $\frac{1}{d_{n,m}}{\frac{{c_{m}}}{m{!}}}S_{n}(\cdot
,H_{m})$ tend to the same limit. From the works of Taqqu \cite{taq79} or
Dobrushin and Major \cite{dob79}, $\frac{1}{d_{n,m}}S_{n}(\cdot ,H_{m})$
converges in law to $Y_{m}(\cdot )$, and so does $\frac{1}{d_{n,m}}%
S_{n}(\cdot ,G-\mathbb{E}\left[ G\left( X\right) \right] )$.\endproof

\end{document}